\documentclass[a4paper,11pt]{article}

\usepackage[utf8]{inputenc}
\usepackage[T1]{fontenc}
\usepackage[english]{babel}
\usepackage{amsmath,amssymb,amsthm,mathtools}
\usepackage[nomath]{lmodern}
\usepackage[margin=1in]{geometry}
\usepackage[babel]{microtype}
\usepackage[pdfusetitle,colorlinks=true,linkcolor=blue,citecolor=blue,urlcolor=blue]{hyperref}
\usepackage[capitalise,noabbrev]{cleveref}

\newtheorem{question}{Question}
\newtheorem{theorem}[question]{Theorem}
\newtheorem{proposition}[question]{Proposition}
\newtheorem{corollary}[question]{Corollary}
\newtheorem{remark}[question]{Remark}

\newcommand{\Z}{\mathbb Z}
\newcommand{\dens}{\operatorname{dens}}
\newcommand{\lcmop}{\operatorname{lcm}}
\newcommand{\abs}[1]{\lvert#1\rvert}

\allowdisplaybreaks
\author{Stijn Cambie\thanks{Department of Computer Science, KU Leuven Campus Kulak--Kortrijk,
8500 Kortrijk, Belgium. Supported by a postdoctoral fellowship of the Research Foundation
Flanders (FWO), grant number 1225224N.}}

\title{Proving it is impossible; on Erd\H{o}s problem \#278}

\begin{document}
\maketitle

\begin{abstract}
Erd\H{o}s and Graham asked for the minimum density missed by one chosen
residue class for each of a prescribed collection of moduli.  We give exact
expressions for natural structured families and relate one of them to hard
balanced-partition instances.  For binary-encoded lists of moduli, allowing
repetitions, even deciding whether the minimum is zero is NP-hard.
\end{abstract}

\section{Introduction}

We consider a question from Erd\H{o}s and Graham
\cite[p.~28, lines 4--6]{ErGr80}, now listed as
\href{https://www.erdosproblems.com/278}{Erd\H{o}s problem \#278}.

\begin{question}[Erd\H{o}s--Graham]
For a finite set of moduli $n_1,\ldots,n_r$, what is the minimum density of
the integers not hit by a suitable choice of congruences
$a_i\pmod {n_i}$?  Is the worst choice obtained by taking all the $a_i$
equal?
\end{question}

The second question asks for the maximum rather than the minimum uncovered
density.  Simpson proved that this maximum is indeed attained when all
residues are equal~\cite[Lemma~2.3]{Simpson86}.  We study the first question.

For a finite list $\mathbf n=(n_1,\ldots,n_r)$, write
\begin{equation}\label{eq:mu}
 \mu(\mathbf n)=
 \min_{0\leq a_i<n_i}
 \dens\!\left(\Z\setminus\bigcup_{i=1}^r(a_i+n_i\Z)\right).
\end{equation}
If $L=\lcmop(n_1,\ldots,n_r)$, every set in \eqref{eq:mu} is periodic
modulo $L$.  Its density is therefore the proportion of its representatives
in $\{0,\ldots,L-1\}$, so the minimum is always well defined.  A tautological
answer is obtained by minimizing over all $\prod_i n_i$ choices of residues;
the point of the problem is to find a more informative description.

We give two complementary reasons why a general description should not be
expected to be simple.  First, the distinct-modulus family
\[
 \{3\}\cup\{3p:p\in P\}
\]
has an exact formula governed by a balanced partition of the primes.  The
same argument gives different multipartition formulas for other common-divisor
patterns.  Moreover, the associated subset-sum order can reproduce the order
type of Chv\'atal's hard knapsack instances~\cite{Chvatal80}; this obstructs a broad class
of natural recursive, comparison-based approaches.

Second, there is a standard complexity-theoretic obstruction.  For the
natural list version, in which a modulus may occur more than once, deciding
whether $\mu(\mathbf n)=0$ is NP-hard.  The proof identifies this question with
Exact Pinwheel Covering, recently proved NP-hard by Kawamura, Kobayashi and
Kusano~\cite{KKK25}.  Hence an exact polynomial-time evaluation of
\eqref{eq:mu} would imply $\mathrm P=\mathrm{NP}$.  This theorem permits
repeated moduli; it does not settle the computational complexity of the
literal pairwise-distinct version.  We keep this distinction explicit below.

\section{A family governed by balanced partition}

\begin{proposition}\label{prop:three}
Let $P$ be a finite set of primes different from $3$.  For the pairwise
distinct moduli
\[
 \mathbf n=\bigl(3,(3p)_{p\in P}\bigr)
\]
one has
\begin{equation}\label{eq:three}
 \mu(\mathbf n)=\frac13
 \min_{P=P_1\mathbin{\dot\cup}P_2}
 \left\{
   \prod_{p\in P_1}\left(1-\frac1p\right)
   +
   \prod_{p\in P_2}\left(1-\frac1p\right)
 \right\}.
\end{equation}
\end{proposition}

\begin{proof}
Translate all congruences so that the chosen class modulo $3$ is
$0\pmod 3$.  A class modulo $3p$ lying in the same fiber adds no newly
covered integer, so in an optimum every such class lies in either
$1\pmod 3$ or $2\pmod 3$.  Let $P_1$ and $P_2$ record these two choices.

Inside the $j$th fiber, the class belonging to $p\in P_j$ excludes one of the
$p$ possible residues.  The Chinese remainder theorem therefore gives the
uncovered proportion
\(
 \prod_{p\in P_j}(1-1/p).
\)
Each fiber has density $1/3$, which proves \eqref{eq:three}.  Conversely,
every partition $P=P_1\mathbin{\dot\cup}P_2$ can be realized by a choice of
congruences.
\end{proof}

Set
\[
 w_p=-\log\left(1-\frac1p\right),
 \qquad W=\sum_{p\in P}w_p.
\]
If $s=\sum_{p\in P_1}w_p$, then the two products in \eqref{eq:three} are
$e^{-s}$ and $e^{-(W-s)}$.  Their sum is minimized precisely when $s$ is as
close as possible to $W/2$.  Thus \eqref{eq:three} is a balanced-partition
problem on the special weights $w_p$.

The same argument also explains why the arithmetic shape of the moduli leads
to different formulas.  If $q$ is prime, $P$ is a set of primes different
from $q$, and the moduli are $(qp)_{p\in P}$, then
\begin{equation}\label{eq:q}
 \mu\bigl((qp)_{p\in P}\bigr)=
 \frac1q
 \min_{P=P_1\mathbin{\dot\cup}\cdots\mathbin{\dot\cup}P_q}
 \sum_{j=1}^q\prod_{p\in P_j}\left(1-\frac1p\right).
\end{equation}
If the modulus $q$ itself is included, one fiber is covered completely and
the analogous minimum is over $q-1$ parts.  Even within \eqref{eq:three},
changing a single prime may change the optimal partition.  There is therefore
no reason for one fixed algebraic expression to describe all common-divisor
patterns.

\section{A precise NP-hardness theorem}\label{sec:NPhard}

Consider the following decision problem, with every integer encoded in binary.

\medskip
\noindent
\textsc{Zero Uncovered Density}

\begin{description}
 \item[Input:] A finite list $\mathbf n=(n_1,\ldots,n_r)$ with $n_i\geq2$.
 \item[Question:] Is $\mu(\mathbf n)=0$?
\end{description}

\begin{theorem}\label{thm:np}
\textnormal{\textsc{Zero Uncovered Density}} is NP-hard under polynomial-time many-one
reductions.
\end{theorem}

\begin{proof}
Exact Pinwheel Covering asks, for a given sequence of positive integers
$n_1,\ldots,n_r$, whether residues $a_i\pmod {n_i}$ can be chosen so that
\begin{equation}\label{eq:cover}
 \bigcup_{i=1}^r(a_i+n_i\Z)=\Z.
\end{equation}
Kawamura, Kobayashi and Kusano proved that this problem is
NP-hard~\cite[Theorem~3]{KKK25}.

It remains only to observe that \eqref{eq:cover} is equivalent to
$\mu(\mathbf n)=0$.  Indeed, fix the residues and let
$L=\lcmop(n_1,\ldots,n_r)$.  If some integer $x$ is uncovered, then every
integer in $x+L\Z$ is uncovered, and hence the uncovered density is at least
$1/L$.  The density is zero exactly when nothing is uncovered.  Thus the
identity map is the desired reduction.

The restriction $n_i\geq2$ causes no problem: an Exact Pinwheel Covering
instance containing the modulus $1$ is a trivial yes-instance and can be
replaced by the fixed yes-instance $(2,2)$.
\end{proof}

\begin{corollary}\label{cor:threshold}
Given a list $\mathbf n$ and a rational $\eta\in[0,1]$, it is NP-hard to
decide whether $\mu(\mathbf n)\leq\eta$, even for the fixed value $\eta=0$.
Consequently, an exact polynomial-time algorithm for \eqref{eq:mu} would imply
$\mathrm P=\mathrm{NP}$.
\end{corollary}

The theorem is deliberately stated as NP-hardness rather than NP-completeness.
In fact, even when a tuple of residues is supplied, deciding whether the
resulting classes cover $\Z$ is coNP-complete
\cite[Theorem~2]{KKK25}.  Thus a proposed solution is not known to have a
polynomial-time verification in the usual sense.

\begin{remark}[The distinct-modulus restriction]\label{rem:distinct}
Exact Pinwheel Covering takes a sequence, and repetitions are allowed.  They
cannot simply be deleted: two copies of the modulus $2$ can choose opposite
residues and cover all integers, whereas one copy cannot.  Therefore
\cref{thm:np} proves hardness for the natural list or multiset extension of
the Erd\H{o}s--Graham problem, but no pairwise-distinct hardness theorem follows
from this reduction.
\end{remark}

\section{Hard balanced partitions with distinct moduli}

We return to the pairwise-distinct family in \cref{prop:three}.  The following
construction retains the original connection with Chv\'atal's hard knapsack
instances while making its scope precise.

For each $n$, the prime number theorem lets us choose distinct primes
$p_1,\ldots,p_n$ in $[X,5X/4]$ once $X$ is large enough.  Since
$-\log(1-1/p)\sim1/p$, we may arrange that, after reordering,
\[
 w_1\geq w_2\geq\cdots\geq w_n>0,
 \qquad
 w_i=-\log\left(1-\frac1{p_i}\right),
 \qquad
 \frac{w_1}{w_n}<\frac32.
\]
Distinct subsets have distinct weight sums.  To see this, equality of two sums would give
\(
 \prod_{i\in I}(1-1/p_i)=\prod_{j\in J}(1-1/p_j).
\)
After canceling $I\cap J$, take the largest prime in $I\cup J$.  On clearing
denominators, that prime divides one side but none of the factors on the other,
a contradiction.

It follows that
\[
 \delta=\min_{I\neq J}
 \abs{\sum_{i\in I}w_i-\sum_{j\in J}w_j}>0.
\]
Choose $x>4n/\delta$, put $m=\lfloor n/2\rfloor$, and define
\[
 c=\lcmop\{\lfloor xw_i\rfloor:1\leq i\leq m\},
 \qquad
 A_i=
 \begin{cases}
  c\lfloor xw_i\rfloor,&i\leq m,\\
  c\lfloor xw_i\rfloor+1,&i>m.
 \end{cases}
\]
Since $\delta\leq w_i$, we have $xw_i>4n$ for every $i$.

\begin{proposition}\label{prop:chvatal}
The weights $(w_i)$ and the integers $(A_i)$ induce the same strict ordering
of all subset sums.  Moreover, the integers $(A_i)$ satisfy the following
properties:
\begin{enumerate}
 \item $\sum_{i\in I}A_i<\frac12\sum_{i=1}^n A_i$ whenever
       $10\abs I\leq n$;
 \item no integer greater than $1$ divides more than $\lceil n/2\rceil$ of
       the $A_i$;
 \item distinct subsets have distinct sums;
 \item no subset has sum $\left\lfloor\frac12\sum_{i=1}^nA_i\right\rfloor$.
\end{enumerate}
\end{proposition}

\begin{proof}
Let $I\neq J$ and cancel their intersection.  With
$\sigma_i=1$ on $I$, $\sigma_i=-1$ on $J$, and $\sigma_i=0$ otherwise, put
$D=\sum_i\sigma_iw_i$.  Then $\abs{xD}>4n$, while
\[
 \sum_i\sigma_iA_i
 =c\sum_i\sigma_i\lfloor xw_i\rfloor
   +\sum_{i>m}\sigma_i.
\]
The difference between the right-hand side and $cxD$ has absolute value less
than $cn+n$.  Hence it has the same sign as $D$, and its absolute value is
greater than $2n$.  This proves the asserted order equivalence and property~3.

Put $y=xw_n>4n$.  The choice of the primes gives
\[
 \frac{\max_i A_i}{\min_i A_i}
 \leq\frac{cxw_1+1}{c(xw_n-1)}
 <\frac{3y/2+1}{y-1}<3.
\]
If $10\abs I\leq n$, then
\[
 \sum_{i\in I}A_i
 \leq\abs I\max_iA_i
 <3\abs I\min_iA_i
 <(n-\abs I)\min_iA_i
 \leq\sum_{i\notin I}A_i,
\]
which is property~1.

For property~2 it is enough to consider a prime divisor $q$.  If $q$ divides
some $A_i$ with $i\leq m$, then, since
$\lfloor xw_i\rfloor$ divides $c$, every prime factor of $A_i$ divides $c$.
But $A_j\equiv1\pmod c$ for $j>m$, so $q$ divides none of the latter
integers.  A prime dividing only integers in the second group divides at most
$\lceil n/2\rceil$ of them as well.

Finally, if the total sum is even, property~4 follows from property~3 by
comparing a subset with its complement.  If the total is odd, a subset with
the stated sum would differ from its complement by $1$, contradicting the
lower bound $2n$ above.
\end{proof}

These are the four structural hypotheses used in Chv\'atal's lower bound for
his broad class of recursive knapsack algorithms~\cite[Theorem~1]{Chvatal80}.
Since the strict order of every pair of subset sums is preserved, comparison
with the complementary subset shows that the feasible subsets, and hence the
optimum, are the same for the two half-capacity knapsack instances.  Thus the
lower bound transfers to the comparison-based members of Chv\'atal's class:
on the prime weights they see exactly the same subset-sum order.  It is
complementary to
\cref{thm:np}: it concerns pairwise-distinct moduli but a restricted algorithmic
model, whereas \cref{thm:np} is standard NP-hardness but permits repetitions.

This qualification is important.  The resemblance of \eqref{eq:three} to
knapsack is not, by itself, a polynomial reduction from an arbitrary knapsack
instance to a set of primes.  The formal many-one reduction in this paper is
the pinwheel reduction of \cref{thm:np}.

\section{Conclusion}

The minimum in Erd\H{o}s problem \#278 can be written explicitly for several
natural common-divisor patterns, but the resulting formulas already contain
balanced multipartition problems.  More decisively, for binary-encoded lists
of moduli, exact evaluation is NP-hard: this holds already at the zero-density
threshold.  Thus no polynomial-time exact answer should be expected unless
$\mathrm P=\mathrm{NP}$.

This does not literally prove that no ``closed formula'' exists, since such a
formula could contain a difficult optimization such as \eqref{eq:three}.
Nor does the standard NP-hardness result presently cover the pairwise-distinct
restriction.  For distinct moduli, \cref{prop:chvatal} instead gives a precise
obstruction for a natural class of recursive methods.  The word
``impossible'' in the title should be understood in these computational senses.

\section*{Acknowledgement}

The author thanks David Pisinger for pointing to the paper~\cite{Chvatal80}. The exact observation on NP-hardness in~\cref{sec:NPhard} is discovered by GPT Pro Sol 5.6, leading to this v2 version.


\begin{thebibliography}{9}

\bibitem{Chvatal80}
V.~Chv\'atal.
\newblock Hard knapsack problems.
\newblock \emph{Operations Research}, 28:1402--1411, 1980.
\newblock \url{https://doi.org/10.1287/opre.28.6.1402}.

\bibitem{ErGr80}
P.~Erd{\H{o}}s and R.~L. Graham.
\newblock \emph{Old and New Problems and Results in Combinatorial Number
Theory}, volume~28 of \emph{Monographies de L'Enseignement Math\'ematique}.
\newblock Universit\'e de Gen\`eve, L'Enseignement Math\'ematique, Geneva,
1980.

\bibitem{KKK25}
A.~Kawamura, Y.~Kobayashi, and Y.~Kusano.
\newblock Pinwheel covering.
\newblock In I.~Finocchi and L.~Georgiadis, editors, \emph{Algorithms and
Complexity---14th International Conference, CIAC 2025, Part II}, volume 15680
of \emph{Lecture Notes in Computer Science}, pages 185--199. Springer, 2025.
\newblock \url{https://doi.org/10.1007/978-3-031-92935-9_12}.

\bibitem{Simpson86}
R.~J. Simpson.
\newblock Exact coverings of the integers by arithmetic progressions.
\newblock \emph{Discrete Mathematics}, 59:181--190, 1986.

\end{thebibliography}
\end{document}